%% file: note.tex
\documentclass[12pt]{article}
\usepackage[T1]{fontenc}
\usepackage{a4wide}

\title{Unified study of the phase transition for block-weighted random planar maps}

\author{Zéphyr Salvy\thanks{Univ Gustave Eiffel, CNRS, LIGM, F-77454 Marne-la-Vallée, France. E-mail: {\tt zephyr.salvy@univ-eiffel.fr}.}}

\makeatletter

\date{}

\usepackage{bbold}
\usepackage{mathtools}
\usepackage{graphicx}
\usepackage{amsmath}
\usepackage{amssymb}
\usepackage{stmaryrd}
\usepackage{comment}
\usepackage{changepage}
\allowdisplaybreaks %to allow page breaks inside equations
\usepackage{microtype} %%% Better typography
\usepackage{zephyr_style_theorems}
\input{macros}

\begin{document}
\thispagestyle{empty}
\maketitle

\begin{abstract}
In [Fleurat, Salvy 2024], we introduced a model of block-weighted random maps that undergoes a phase transition as the density of separating elements changes. The purpose of this note is to demonstrate that the methodology we developed can be extended to many other families of maps. We prove that a phase transition exists and provide detailed information about the size of the largest blocks in each regime.

{\noindent\small\bf DOI:}  {\tt https://doi.org/10.5817/CZ.MUNI.EUROCOMB23-109}
\end{abstract}

\maketitle
\input{note_body}

\bibliographystyle{alpha}
\bibliography{note.bib}

\section*{Acknowledgements}
The author would like to thank their supervisors Marie Albenque and Éric Fusy for their guidance and support.

%\section*{Affiliation}
\end{document}

%% file: macros.tex
\renewcommand{\leq}{\leqslant}
\renewcommand{\geq}{\geqslant}
\renewcommand{\(}{\left(}
\renewcommand{\)}{\right)}

\renewcommand{\epsilon}{\varepsilon}

\newcommand{\N}{\mathbb{Z}_{\geq0}}

\newcommand{\m}{\mathfrak{m}}

\newcommand{\M}{\mathbf{M}}

\newcommand{\pr}[2][]{\mathbb{P}#1\( #2\)}
\newcommand{\E}[1]{\mathbb{E}\left[ #1 \right]}

%% file: note_body.tex
% !TEX root = note.tex
\section{Introduction}
A \emph{planar map} $\m$ is the proper embedding into the two-dimensional sphere of a connected planar finite multigraph, considered up to homeomorphisms. Maps exhibit very rich combinatorial and probabilistic properties, which have been the focus of an extensive literature. Many families of planar maps have very nice counting formulas \cite{tutte_1963}. A key aspect of planar maps is that they can be decomposed, typically into components of higher connectivity degree. Such decompositions typically relate one family of planar maps to another and gives an equation between their generating series.

Theses types of decompositions were initially introduced by Tutte \cite{tutte_1963} to obtain some enumerative results about planar maps. But they also play a major role in the enumerative study of planar graphs \cite{GN09}. They allow to study certain models of discrete metric spaces in theoretical physics \cite{bonzomLagrange}. In view of applications to random generation \cite{random_sampling_maps}, the decomposition of planar maps has been systematised in \cite{airy}, where a uniform treatment via analytic combinatorics is developed. A probabilistic approach was later derived using encoding of maps via enriched trees \cite{Stu18,2Louigi}.

Planar map models exhibit \emph{universality}, meaning that many natural classes of random maps show similar behavior as their size grows to infinity. When taking an object of size $n$ uniformly at random among all objects in a class and appropriately rescaling its distance, the sequence of random objects converges to a certain metric space. This was first proved for uniform quadrangulations by Miermont \cite{Miermont2013} and Le Gall \cite{LeGall2013}, and since then, results have been extended to other families of maps, including uniform triangulations and uniform $2q$-angulations ($q\geq2$) \cite{LeGall2013}, uniform simple triangulations and uniform simple quadrangulations \cite{add-alb}, bipartite planar maps with a prescribed face-degree sequence \cite{Marzouk2018}, $(2q+1)$-angulations \cite{AddarioBerryAlbenque2021} and Eulerian triangulations \cite{Carrance2021}.

In a previous article \cite{FleuratSalvy23}, together with Fleurat, we studied a model of random maps, depending on a parameter $u$ which controls the density of separating elements. We proved that this model exhibits a phase transition as $u$ varies, and that it interpolates between the Brownian sphere and the Brownian tree of Aldous \cite{CRTII}. This approach~---~which we detailed in \cite{FleuratSalvy23} for general maps and their $2$-connected cores and for general quadrangulations and their simple cores~---~can be applied to other decompositions, such as those in \cite[Table 3]{airy} (which is partially reproduced in \cref{airy-table-3}), and this is the focus of this note. We restrict our study to decomposition schemes without ``coreless'' maps (the decompositions involving coreless maps, such as $2$-connected maps into $3$-connected components, bring further difficulties, which we expect to handle with some more work).

Let us give some formalism for decompositions. A map is said to be \emph{loopless} if it does not contain any loop; \emph{$2$-connected} if it does not contain any \emph{cut vertex} (\emph{i.e.} a vertex whose removal deconnects the map) and \emph{simple} if it has neither loops nor multiple edges. Planar maps can be decomposed into loopless (or $2$-connected, or simple, or $2$-connected simple...) components, which are the so-called ``blocks''. It is also the case for \emph{bipartite} maps, whose vertices can be properly bicolored in black and white; and for \emph{triangulations}, whose faces all have degree $3$. The latter can be decomposed into \emph{irreducible} components, in which every $3$-cycle defines a face. We consider eight models here (see \cref{airy-table-3}):
\begin{enumerate}
\item Loopless maps decomposed into simple blocks;
\item General maps decomposed into $2$-connected blocks;
\item $2$-connected maps decomposed into $2$-connected simple blocks;
\item Bipartite maps decomposed into bipartite simple blocks;
\item Bipartite maps decomposed into bipartite $2$-connected blocks;
\item Bipartite $2$-connected maps decomposed into bipartite $2$-connected simple blocks;
\item Loopless triangulations decomposed into triangular simple blocks;
\item Simple triangulations decomposed into triangular irreducible blocks.
\end{enumerate}

In general, the \emph{size} $|\m|$ of a planar map $\m$ is its number of edges.
In a decomposition scheme, we let $M(z) = \sum_{n\in\N} m_n z^n$
be the generating series of the class of maps to be decomposed and similarly
$B(z) = \sum_{n\in\N} b_n z^n$
be the generating series of the class of ``blocks'' into which the maps are decomposed.

Setting $M(z,u)=\sum_{\m\in \mathcal{M}}z^{|\m|}u^{b(\m)}$, where $b(\m)$ is the number of blocks of positive size in $\m$, all the models we consider~---~which are listed in \cref{airy-table-3}~---~satisfy
\begin{equation}
\label{equation-decomp-1}
M(z,u) = uB(H(z, M(z,u))),
\end{equation}
or, for the last one,
\begin{equation}
\label{equation-decomp-2}
M(z,u) = (1+M(z,u)) \times u B(H(z, M(z,u))).
\end{equation}
For example, for the decomposition of general maps into $2$-connected ones (which is the case studied in \cite{FleuratSalvy23}), one has
\[M(z,u) = u B(z(1+M(z,u))^2).\]
For $u>0$, denote by $\rho(u)$ the radius of convergence of $z\mapsto M(z,u)$. In view of the form of \cref{equation-decomp-1,equation-decomp-2} and in particular that they are non-linear, it holds that $M(\rho(u), u)<\infty$.

In the following, $\M$ and $\M_n$ are random variables drawn according to the following probability distributions. For $u\in \mathbb{R}_{>0}$, $n\in \mathbb{Z}_{>0}$ and $\m \in \mathcal{M}$, we set
\[\pr[_{u}]{\m} = \frac{\rho(u)^{|\m|}u^{b(\m)}}{M(\rho(u),u)}\qquad\text{and}\qquad\pr[_{n,u}]{\m} = \frac{u^{b(\m)}}{[z^n] M(z,u)}\mathbb{1}_{|\m|=n}.\]

\begin{table}
\begin{center}
\begin{tabular}{c|llc}
Scheme & maps, $M(z)$ & blocks, $B(z)$ & submaps, $H(z, M)$\\
\hline
1 & loopless, $M_2(z)$ & simple, $M_3(z)$ & $z(1+M)$ \\
2 & all, $M_1(z)$ & $2$-connected, $M_4(z)$ & $z(1+M)^2$\\
3 & $2$-connected $M_4(z)-z$ & $2$-connected simple, $M_5(z)$ & $z(1+M)$\\
\hline
4 & bipartite, $B_1(z)$ & bipartite simple, $B_2(z)$ & $z(1+M)$\\
5 & bipartite, $B_1(z)$ & bipartite $2$-connected, $B_4(z)$ & $z(1+M)^2$\\
6 & bipartite $2$-connected, $B_4(z)$ & bipartite $2$-connected simple $B_5(z)$ & $z(1+M)$\\
\hline
7 & loopless triangulations, $T_1(z)$ & simple triangulations, $z+zT_2(z)$ & $z(1+M)^3$\\
8 & simple triangulations, $T_2(z)$ & irreducible triangulations, $T_3(z)$ & $z(1+M)^2$\\
\end{tabular}
\vspace{\abovecaptionskip}
\caption{Partial reproduction of \cite[Table 3]{airy}, which describes composition schemas of the form $\mathcal{M} = \mathcal{B}\circ\mathcal{H}$ except the last one where $\mathcal{M} = (1+\mathcal{M})\times \mathcal{B}\circ\mathcal{H}$. The terminology and notation were slightly changed. For all $i$, $[z^n]M_i(z)$ and $[z^n]B_i(z)$ is the number of such maps with $n$ edges. $[z^n]T_1(z)$ (resp. $[z^n]T_2(z)$ and $[z^n]T_3(z)$) is the number of loopless (resp. simple or irreducible) triangulations with $n+2$ (resp. $n+3$) vertices.}
\label{airy-table-3}
\end{center}
\end{table}

Regarding enumeration in our setting, we show the following by analytic methods (details are omitted in this short note).

\begin{theorem}
\label{th-enum-models-airy}
For any model described in \cref{airy-table-3}, where maps are decomposed into blocks weighted with a weight $u>0$, there exists a critical value $u_C$ at which the model undergoes a phase transition. As $u$ varies, there exists $c(u) > 0$ such that
    \[[z^n] M(z, u) \sim \begin{cases}c(u) n^{-5/2} \rho(u)^{-n}\quad\text{if}\ u < u_C\\
    c(u_C) n^{-5/3} \rho(u_C)^{-n}\quad\text{if}\ u = u_C\\
    c(u) n^{-3/2} \rho(u)^{-n}\quad\text{if}\ u > u_C
    \end{cases}.\]
\end{theorem}
All the constants involved in \cref{th-enum-models-airy} are explicit. \Cref{values-models-airy} gives the expressions for $u_C$, $\rho(u)$ and $M(\rho(u),u)$ when $u\leq u_C$.

The polynomial correction for $u<u_{C}$ (subcritical case) is the same than for planar maps, whereas when $u>u_C$ (supercritical case) it is the same than for plane trees. Moreover, at $u=u_{C}$, a new asymptotic behaviour appears with a polynomial correction in $n^{-5/3}$.

In this note, we also focus on another aspect of the phase transition, namely the size of the largest blocks. We show that if $u < u_C$, a condensation phenomenon occurs and the largest block is of size $\Theta(n)$; when $u > u_C$, the largest block is of size $\Theta(\log(n))$; for $u = u_C$, the largest block is of size $\Theta(n^{2/3})$ (\cref{size-largest-block}). For the subcritical case, as in \cite{FleuratSalvy23}, we follow the probabilistic approach of \cite{2Louigi} (whereas \cite{airy} gives an analytic approach).

These results further support that the scaling limits should be the Brownian sphere when $u<u_C$, the Brownian tree when $u>u_C$ and the stable tree of parameter $3/2$ when $u=u_C$. This was proved for the decomposition of quadrangulations into simple components \cite{FleuratSalvy23}, and we expect this phenomenon to be generic. For model 2, the critical scaling limit was established in \cite{FleuratSalvy23} and the supercritcal one in \cite{Stufler-survey-2020}. For model 5, the supercritical case was also established \cite{Stufler-survey-2020}.

\section{Tree structure}
We explain here how an underlying tree structure can be associated to each of the models of \cref{airy-table-3}. As a first step, we rewrite the decomposition equations in the standard Lagrangian form $M(z) = z \times \Phi(M(z))$ for some function $\Phi$, taking the weight $u$ into account. (Beware that \cref{equation-decomp-1,equation-decomp-2} are not of this form as the products by $z$ are inside $H$.)

\begin{proposition}\label{prop-lagrangian}
For all models listed in \cref{airy-table-3}, there exists a generating function $\Phi$ with nonnegative coefficients such that
\begin{equation}
\label{lagrangian-tree}
M(z,u) = z \times \Phi(M(z,u),u).
\end{equation}
\end{proposition}

\begin{figure}
\begin{center}
\includegraphics[height=0.25\textheight, width=\textwidth]{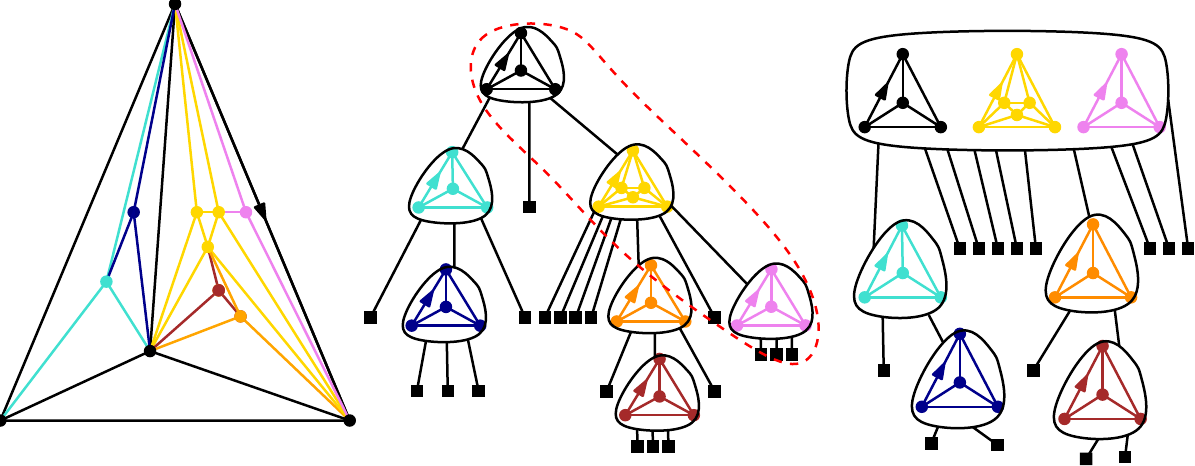}
\end{center}
\caption{A simple triangulation, its classical tree of (irreducible) blocks, the adapted tree where some blocks are grouped into sequences.}
\label{tree-decomp-simple-irred-tri}
\end{figure}

\begin{proofsketch}
\label{def-d}
We discuss how the rewriting is done for two cases: general maps into $2$-connected components, and simple triangulations into irreducible components.

We start by maps decomposed into $2$-connected components. In order to do the rewriting, we need to change the size parameter, which we take as the number of half-edges plus one. Accordingly, we set $\hat{M}(z,u) = z(1+ M(z^2,u))$. The equation then becomes
\[\hat{M}(z,u) = z \(1 + u B(\hat{M}(z,u)^2)\),\]
which is of the desired form.

The second case we discuss is the decomposition of simple triangulations into irreducible components. We also need to change the size parameter, taking here the number of inner faces instead of the number of vertices. A further ingredient compared to the first case is that we need to group the components into sequences to obtain an equation in Lagrangian form.
Let $\hat{z}$ count internal faces. The generating series of simple (resp. irreducible) triangulations counted by internal faces $\hat{T_2}(\hat{z},u)$ (resp. $\widetilde{T_3}(\hat{z})$) is closely related to $T_2(z,u)$ (resp. $T_3(z)$) since a triangulation with $n+3$ vertices has $2n+2$ faces so $2n+1$ internal faces. Then, denoting $\hat{T_3}(\hat{z}) = \widetilde{T_3}(\hat{z})/\hat{z}$, It holds that
\[\hat{T_2}(\hat{z},u)= \hat{z} + u \hat{T_2}(\hat{z},u) \hat{T_3}(\hat{T_2}(\hat{z},u)),\qquad\text{so}\qquad\hat{T_2}(\hat{z},u) = \frac{\hat{z}}{1 - u \hat{T_3}(\hat{T_2}(\hat{z},u))}.\]
Therefore, we set $\Phi(M,u) = \frac{1}{1 - u \hat{T_3}(M)}$. This corresponds to the vertices of the tree encoding a sequence of irreducible triangulations, which is represented on \cref{tree-decomp-simple-irred-tri}\footnote{With additional work, one can do the same for simple quadrangulations decomposed into irreducible ones.}.
\end{proofsketch}

The \emph{block tree} $T_{\m}$ of a map $\m$ is the tree associated to the the decomposition of $\m$ as expressed in \cref{prop-lagrangian}. Each node of the tree correspond to an object $\varphi$ counted by $\Phi$. The subtrees hanging at a node (corresponding to some $\varphi$) are the trees of the components substituted into $\varphi$.

As a consequence of \cref{lagrangian-tree}, if we set a probability measure $\mu^{u}$ such that, for $k\in\N$,
\[\mu^{u}(k) = \frac{[X^k] \Phi(X, u) y(u)^k}{\Phi(y(u),u)}\]
for $y(u) = M(\rho(u),u)$ (using the definition of $M$ given in \cref{lagrangian-tree}); then, using the fact that $T_\m$ and the decoration of its vertices are bijectively linked to $\m$, we have the following result (again stated in terms of the $M$ of \cref{lagrangian-tree}):

\begin{theorem}
\label{tree-structure}
For all $u>0$, $T_\M$ follows the law of a Galton-Watson tree of reproduction law $\mu^u$.
Moreover, $T_{\M_n}$ follows the law of a Galton-Watson tree of reproduction law $\mu^u$ conditioned to have $n$ vertices.
\end{theorem}

The decomposition tree of a random map of ``Lagrangian size'' $n$ is a Galton-Watson tree of reproduction law $\mu^u$ conditioned to have $n$ vertices. For instance, for maps decomposed into $2$-connected components, the tree of a random map with $n$ edges is a Galton-Watson tree conditioned to have $2n+1$ vertices. This enables to put into light a phase transition on the tree structure, using the usual phase transition for Galton-Watson trees \cite{GW}.

\begin{proposition}
\label{prop:expectation-mu-u}
The expectation $E(u)$ of $\mu^u$ is written down in \cref{values-models-airy}.
\end{proposition}

\begin{table}
\begin{adjustwidth}{-1cm}{}
\begin{center}
\def\arraystretch{1.5}
\begin{tabular}{c|c|cc|cc}
Scheme & $u_C$ & $\rho(u)$ & $M(\rho(u),u)$ & $E(u)$ & $1-E(1)$\\
\hline
1 & $\frac{81}{17}$ & $\frac{27}{8(5u + 27)}$ & $\frac{5u}{27}$ & $\frac{32 u}{3 (5 u+27)}$ & $\frac{2}{3}$\\
2 & $\frac{9}{5}$ & $\frac{4}{3(u^2 + 6u + 9)}$ & $\frac{u}{3}$ & $\frac{8u}{3(u + 3)}$ & $\frac{1}{3}$\\
3 & $\frac{135}{7}$ & $\frac{128}{27(5u + 27)}$ & $\frac{25u^2 + 135u + 128}{27(5u + 27)}$ & $\frac{32u}{5(5u + 27)}$ & $\frac{4}{5}$\\
\hline
4 & $\frac{36}{11}$ & $\frac{5}{8(u + 4)}$ & $\frac{u}{4}$ & $\frac{20u}{9(u + 4)}$ & $\frac{5}{9}$ \\
5 & $\frac{52}{27}$ & $\frac{25}{8(u^2 + 8u + 16)}$ & $\frac{u}{4}$ & $\frac{40u}{13(u + 4)}$ & $\frac{5}{13}$ \\
6 & $\frac{68}{3}$ & $\frac{125}{128(u + 4)}$ & $\frac{u}{4}$ & $\frac{20u}{17(u + 4)}$ & $\frac{13}{17}$\\
\hline
7 & $\frac{16}{7}$ & $\frac{54}{u^3 + 24u^2 + 192u + 512}$ & $\frac{u}{8}$ & $\frac{9u}{2(u + 8)}$ & $\frac{1}{2}$\\
8 & $\frac{64}{37}$ & $\frac{25}{6912}u^2 - \frac{5}{108}u + \frac{4}{27}$ & $\frac{5u}{32 - 5u}$ & $\frac{27u}{2(32 - 5u)}$ & $\frac{1}{2}$\\
\end{tabular}
\vspace{\abovecaptionskip}
\caption{Values of $u_C$, $\rho(u)$, $M(\rho(u), u)$ and $E(u)$ when $u\leq u_C$ for all the decomposition schemes of \cref{airy-table-3}.}
\label{values-models-airy}
\end{center}
\end{adjustwidth}
\end{table}

\section{Results on the size of the largest blocks}
Starting from \cref{equation-decomp-1,equation-decomp-2,tree-structure}, we use the same techniques as in \cite{FleuratSalvy23} to obtain results on decomposition schemes. However, simple triangulations decomposed into irreducible blocks present a challenge, as the vertices of the block trees are not decorated with a single block (or none) but with a sequence of blocks. Hence, the size of the blocks cannot be immediately read from the degrees in the block tree. However, an extreme condensation phenomenon occurs, concentrating mass in only one element of the sequence (as in \cite[Th1]{Gourdon98}), resulting in a similar behaviour.

Denote by $L_{n,j}$ the size of the $j$-th largest block of $\M_n$. By \cref{tree-structure}, the same arguments as in \cite{FleuratSalvy23} apply and the following holds.

\begin{theorem}
\label{size-largest-block}
Models described in \cref{airy-table-3}, where maps are decomposed into blocks weighted with a weight $u>0$, satisfy the following.

\begin{description}
    \item[Subcritical case] For all $u < u_C$, we have
\[L_{n,1} = (1-E(u))n + O_{\mathbb{P}}(n^{2/3})\quad\text{and}\quad L_{n,2}=O_{\mathbb{P}}(n^{2/3}).\]
Moreover, there exists an explicit constant $\tilde{c}(u)>0$ such that the following joint convergence holds:
\begin{equation}
\label{th-sous-crit}
\(\frac{1}{n \tilde{c}(u)}\)^{2/3}\((1-E(u))n-L_{n,1}, \(L_{n,j}, j\geq 2\)\)\xrightarrow[n\to\infty]{(d)} \(L_1,\(\Delta L_{(j-1)},j\geq 2\)\)
\end{equation}
where $(L_t)_{t\in[0,1]}$ is a Stable process of parameter $3/2$ such that $\E{e^{-sL_1}} = e^{\Gamma(-3/2)s^{3/2}}$
and $\Delta L_{(1)} \geq \Delta L_{(2)} \geq \dots$ is the ranked sequence of its jumps.

\item[Supercritical case] For all $u > u_C$, there exist explicit values $F(u), G(u) > 0$ such that, for all fixed $j\geq 1$,
\[L_{n,j} = F(u)\ln(n) - G(u) \ln(\ln(n))+O_{\mathbb{P}}(1).\]

\item[Critical case] If $u = u_C$, then
\[\(\frac{L_{n,j}}{n^{2/3}}, j \geq 1\) \xrightarrow[n\to\infty]{(d)} \(E_{(j)}, j \geq 1\),\]
where the $\(E_{(j)}\)$ are the ordered atoms of an explicite Point Process, specified in \cite[Ex19.27, Rk19.28]{survey-trees}.
\end{description}
\end{theorem}

As mentioned in \cite{FleuratSalvy23}, for $u=1$, we retrieve by a probabilistic method the results of \cite[Table 4]{airy}, established by analytic techniques: indeed, our $1-E(1)$ corresponds to their $\alpha_0$. The probabilistic approach we follow was first developed by \cite{2Louigi} and has the advantage that we obtain a joint limit law for the largest block and the subsequent ones. It can also yield local limit theorems for the size of the largest block, as is discussed by Stufler in \cite{Stu20}.